\documentclass[12pt]{amsart}

\usepackage{enumitem}
\usepackage{psfrag}
\usepackage{graphicx}
\usepackage{pinlabel}
\usepackage{graphicx}
\usepackage{amsmath}
\usepackage{amssymb}
\usepackage{amscd}
\usepackage{autobreak}
\usepackage{picinpar}
\usepackage{times}
\usepackage{pb-diagram}
\usepackage{graphicx}
\usepackage{wrapfig}
\usepackage{xspace}
\usepackage{hyperref}
\usepackage{url}
\usepackage{caption}
\usepackage{subcaption}
\usepackage{fancyhdr}
\usepackage{overpic}
\usepackage{color}
\usepackage[square,comma,sort&compress,numbers]{natbib}
\usepackage{latexsym}
\usepackage{mathrsfs}
\usepackage{amsfonts}
\usepackage{hyperref}
\usepackage{amsthm}
\usepackage{appendix}
\usepackage{pdfsync}

\theoremstyle{definition}
\newtheorem{theorem}{Theorem}[section]
\newtheorem{lemma}[theorem]{Lemma}
\newtheorem{proposition}[theorem]{Proposition}
\newtheorem{corollary}[theorem]{Corollary}
\newtheorem{definition}[theorem]{Definition}

\newtheorem{remark}[theorem]{Remark}
\newtheorem*{theorem*}{Theorem}
\def\qed{\hfill{Q.E.D.}\smallskip}

\begin{document}

\title{\bf Rigidity of generalized Thurston's sphere packings on 3-dimensional manifolds with boundary}
\author{Xu Xu, Chao Zheng}

\date{\today}

\address{School of Mathematics and Statistics, Wuhan University, Wuhan, 430072, P.R. China}
 \email{xuxu2@whu.edu.cn}

\address{School of Mathematics and Statistics, Wuhan University, Wuhan 430072, P.R. China}
\email{czheng@whu.edu.cn}


\begin{abstract}
Motivated by Guo-Luo's generalized circle packings on surfaces with boundary \cite{GL2},
we introduce the generalized Thurston's sphere packings on 3-dimensional manifolds with boundary.
Then we investigate the rigidity of the generalized Thurston's sphere packings.
We prove that the generalized Thurston's sphere packings are locally determined by the combinatorial scalar curvatures.
We further prove the infinitesimal rigidity that the generalized Thurston's sphere packings can not be deformed while keeping the combinatorial Ricci curvatures fixed.
\end{abstract}

\maketitle

\section{Introduction}
This is a continuation of \cite{XZ}, in which we investigated the rigidity and deformation of the generalized tangential sphere packings on 3-dimensional manifolds with boundary.
In this paper, motivated by Guo-Luo's generalized circle packings on surfaces with boundary \cite{GL2},
we introduce the generalized Thurston's sphere packings on 3-dimensional manifolds with boundary.
These are analogues of Thurston's Euclidean sphere packings on 3-dimensional closed manifolds studied in \cite{DG, Glickenstein JDG, HX} and
natural generalizations of the generalized tangential sphere packings on 3-dimensional manifolds with boundary studied in \cite{XZ}.
The main focus of this paper is to study the rigidity of the generalized Thurston's sphere packings on 3-dimensional manifolds with boundary.

Suppose $\widetilde{\Sigma}$ is a compact 3-dimensional manifolds with the boundary $\partial\widetilde{\Sigma}$ consisting of $N$ connected components.
By coning off each boundary component of $\widetilde{\Sigma}$ to be a point,
one can obtain a compact 3-dimensional space, denoted by $\Sigma$.
There are exactly $N$ cone points $\{v_1,...,v_N\}$ in $\Sigma$ and $\Sigma-\{v_1,...,v_N\}$ is homeomorphic to $\widetilde{\Sigma}-\partial \widetilde{\Sigma}$.
An ideal triangulation $\mathcal{T}$ of $\widetilde{\Sigma}$ is a triangulation $\mathcal{T}$ of $\Sigma$ such that the vertices of the triangulation are exactly the cone points $\{v_1,...,v_N\}$.
Then $(\Sigma,\mathcal{T})$ is called a closed 3-dimensional manifold with ideal triangulation.
Replacing each tetrahedron in $\mathcal{T}$ by a hyper-ideal tetrahedron and replacing the affine gluing homeomorphisms by isometries preserving the corresponding hyperbolic hexagonal faces,
we obtain a hyper-ideal polyhedral metric
on $(\Sigma,\mathcal{T})$.
We call $(\Sigma,\mathcal{T})$ an ideally triangulated compact 3-dimensional manifold with boundary.

\begin{figure}[!ht]
\centering
\begin{overpic}[scale=1]{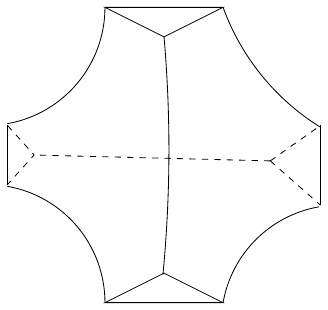}
\put (100,42) {$\triangle_h$}
\put (48,-5) {$\triangle_k$}
\put (-10,45) {$\triangle_j$}
\put (60,96) {$\triangle_i$}
\put (48,97) {$x^i_{jh}$}
\put (56,80) {$x^i_{kh}$}
\put (37,80) {$x^i_{jk}$}
\put (18,72) {$l_{ij}$}
\put (55,60) {$l_{ik}$}
\put (16,20) {$l_{jk}$}
\put (80,72) {$l_{ih}$}
\put (37,40) {$l_{jh}$}
\put (82,19) {$l_{kh}$}
\end{overpic}
\caption{Hyper-ideal tetrahedron $\sigma=\{ijkh\}$}
\label{Figure 2}
\end{figure}
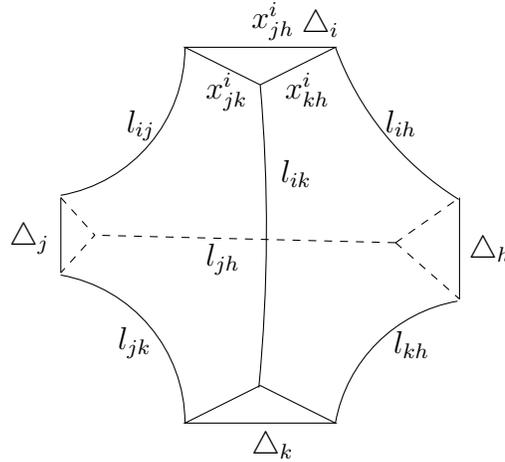

A hyper-ideal tetrahedron is a compact convex polyhedron in $\mathbb{H}^3$ that is diffeomorphic to a truncated tetrahedron in $\mathbb{E}^3$,
which has four right-angled hyperbolic hexagonal faces and four hyperbolic triangular faces.
Any triangular face is required to be orthogonal to its three adjacent hexagonal faces.
The four triangular faces isometric to hyperbolic triangles are called vertex triangles.
An edge in a hyper-ideal tetrahedron is the intersection of two hexagonal faces and a vertex edge is the intersection of a hexagonal face and a vertex triangle.
We use $E$, $F$ and $T$ to represent the set of edges, faces and hyper-ideal tetrahedra in the triangulation $\mathcal{T}$ respectively.
A hyper-ideal tetrahedron with four vertex triangles $\triangle_\nu, \nu=i,j,k,h$ is denoted by $\sigma=\{ijkh\}$.
Using the Klein model of $\mathbb{H}^3$,
the hyper-ideal tetrahedron $\sigma=\{ijkh\}$ corresponds to a Euclidean tetrahedron with four vertices $i,j,k,h$.
We call these vertices as hyper-ideal vertices of the hyper-ideal tetrahedron $\sigma=\{ijkh\}$, which corresponds to vertex triangles $\triangle_\nu, \nu=i,j,k,h$.
The set of all hyper-ideal vertices in the hyper-ideal tetrahedra is denoted by $V$.
The edge joining $\triangle_i$ to $\triangle_j$ is denoted by $\{ij\}\in E$.
The hexagonal face adjacent to $\{ij\}$, $\{jk\}$ and $\{ik\}$ is denoted by $\{ijk\}\in F$.
The length of the edge $\{ij\}\in E$ is denoted by $l_{ij}$.
The length of the vertex edge $\triangle_i\bigcap\{ijk\}$ is denoted by $x^i_{jk}$.
Please refer to Figure \ref{Figure 2}.

\begin{definition}\label{Def: SP metric}
Suppose $(\Sigma,\mathcal{T})$ is an ideally triangulated compact 3-dimensional manifold with boundary.
Let $\Phi: E\rightarrow [0,\frac{\pi}{2}]$ be a weight defined on the edges of $(\Sigma,\mathcal{T})$.
A generalized Thurston's sphere packing metric on $(\Sigma,\mathcal{T},\Phi)$ is defined to be a map $r: V\rightarrow (0,+\infty)$ such that
\begin{description}
\item[(1)] The length $l_{ij}$ of the edge $\{ij\}\in E$ between two vertex triangles $\triangle_i$ and $\triangle_j$ is given by
\begin{equation}\label{Eq: SP metric}
l_{ij}=\cosh^{-1}(\sinh r_i\sinh r_j+\cos\Phi_{ij}\cosh r_i\cosh r_j).
\end{equation}
\item[(2)] The lengths $l_{ij}, l_{ik}, l_{ih}, l_{jk}, l_{jh}, l_{kh}$ determine a hyper-ideal tetrahedron $\{ijkh\}$.
\end{description}
\end{definition}

The function $l: E\rightarrow (0, +\infty)$ determined by (\ref{Eq: SP metric}) is a discrete hyperbolic metric on $(\Sigma,\mathcal{T})$.
The condition (2) in Definition \ref{Def: SP metric} is called the non-degenerate condition.
Intuitively, a degenerate hyper-ideal tetrahedron is a collapsed octagon, which was called a flat hyper-ideal tetrahedron by Luo-Yang \cite{L-Y}.
Furthermore, Luo-Yang \cite{L-Y} showed that a hyper-ideal tetrahedron is non-degenerate if and only if the four vertex triangles are non-degenerate, i.e., triangle inequalities hold for any vertex triangle.
Please see Section \ref{section 2} or \cite{L-Y} for more details.

\begin{figure}[!ht]
\centering
\begin{overpic}[scale=1]{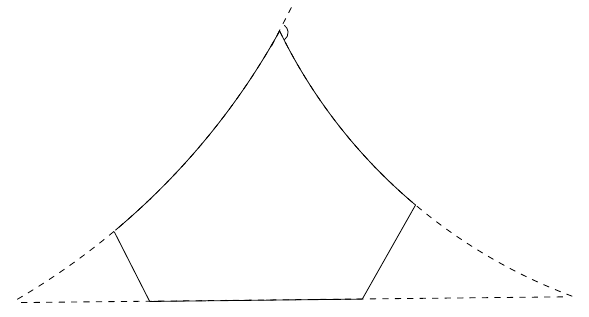}
\put (60,32) {$r_j$}
\put (30,32) {$r_i$}
\put (0,0) {$i$}
\put (98,0) {$j$}
\put (45,-2) {$l_{ij}$}
\put (50,46) {$\Phi_{ij}$}
\put (17,6) {$\theta_{ij}$}
\put (67,7) {$\theta_{ji}$}
\end{overpic}
\caption{Hyperbolic pentagon}
\label{Figure 1}
\end{figure}
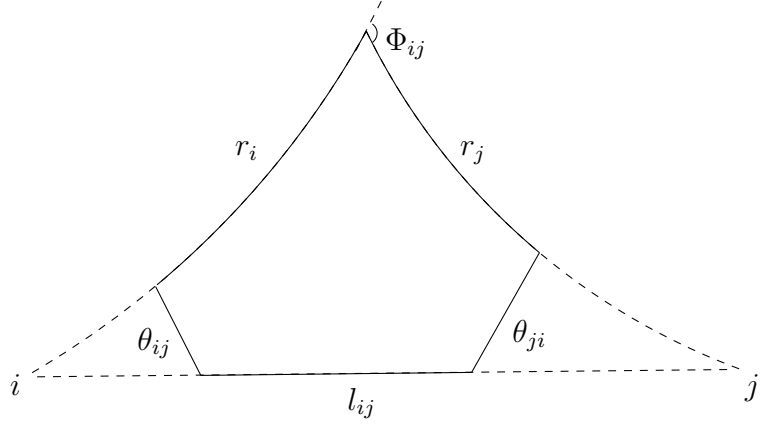

If $\Phi\equiv0$, then the generalized Thurston's sphere packing metric is exactly the generalized tangential sphere packing metric introduced in \cite{XZ}.
We give a geometric explanation of the generalized Thurston's sphere packing metric.
Given a generalized Thurston's sphere packing metric $r$ on $(\Sigma,\mathcal{T},\Phi)$,
a sphere $B_i$ of radius $r_i$ is attached to each hyper-ideal vertex $i\in V$, and the intersection angle of $B_i$ and $B_j$ attached to the hyper-ideal vertices $i,j\in V$ is $\Phi_{ij}$.
Then the length $l_{ij}$ of the edge $\{ij\}$ is a function of $r_i,\ r_j$ and the fixed $\Phi_{ij}$ defined by (\ref{Eq: SP metric}) by the cosine law for hyperbolic pentagons.
Please refer to Figure \ref{Figure 1}.

For each $\sigma\in T$, the generalized Thurston's sphere packing metric determines a non-degenerate hyper-ideal tetrahedron by Definition \ref{Def: SP metric}.
Gluing these hyper-ideal tetrahedra along the faces by isometries may produce singularities.
We use the following combinatorial Ricci curvature to describe the singularities along the edges.

\begin{definition}[\cite{Glickenstein JDG}]\label{Def: CRC}
Suppose $(\Sigma,\mathcal{T})$ is an ideally triangulated compact 3-dimensional manifold with boundary
and $l: E\rightarrow(0,+\infty)$ is a discrete hyperbolic metric on $(\Sigma,\mathcal{T})$.
The combinatorial Ricci curvature of $l$ along the edge $\{ij\}\in E$ is defined to be
\begin{equation}\label{Eq: CRC}
K_{ij}=2\pi-\sum_{\{ijkh\}\in T}\beta_{ij,kh},
\end{equation}
where $\beta_{ij,kh}$ is the dihedral angle of the hyper-ideal tetrahedron $\{ijkh\}$ at the edge $\{ij\}\in E$ and the summation is taken over the hyper-ideal tetrahedra with $\{ij\}$ as a common edge.
\end{definition}

We have the following infinitesimal rigidity for the generalized Thurston's sphere packing metrics.
\begin{theorem}\label{Thm: infinitesimal rigidity}
Suppose $(\Sigma,\mathcal{T})$ is an ideally triangulated compact 3-dimensional manifold with boundary.
Let $\Phi: E\rightarrow [0,\mathrm{arccos}\frac{1}{3}]$ be a weight defined on the edges of $(\Sigma,\mathcal{T})$.
Then the generalized Thurston's sphere packing metric on $(\Sigma,\mathcal{T},\Phi)$ can not be deformed while keeping the combinatorial Ricci curvatures along the edges fixed.
\end{theorem}

Similar to the combinatorial scalar curvature for Thurston's Euclidean sphere packing metrics on 3-dimensional closed manifolds \cite{Glickenstein JDG, HX}, we introduce the following combinatorial scalar curvature for the generalized Thurston's sphere packing metrics
to describe the singularities at the hyper-ideal vertices.

\begin{definition}\label{Def: CSC}
Suppose $(\Sigma,\mathcal{T})$ is an ideally triangulated compact 3-dimensional manifold with boundary.
Let $\Phi: E\rightarrow [0,\frac{\pi}{2}]$ be a weight defined on the edges of $(\Sigma,\mathcal{T})$.
The combinatorial scalar curvature at a hyper-ideal vertex $i\in V$ for a generalized Thurston's sphere packing metric $r$ is defined to be
\begin{equation*}
K_i=\sum_{j;j\sim i}K_{ij}\cosh \theta_{ij},
\end{equation*}
where $\theta_{ij}$ is the edge length opposite to the edge with length $r_j$ in the hyperbolic pentagon formed by three edges with lengths $r_i,\ r_j,\ l_{ij}$ as shown in Figure \ref{Figure 1}.
\end{definition}

\begin{remark}
If $\Phi_{ij}=0$ for any edge $\{ij\}\in E$, then $\cosh \theta_{ij}=1$ and hence
\begin{equation*}
K_i=\sum_{j;j\sim i}\bigg(2\pi-\sum_{\{ijkh\}\in T}\beta_{ij,kh}\bigg)
=2\pi\chi(\Sigma_i)+\mathrm{Area}(\Sigma_i),
\end{equation*}
where $\Sigma_i$ is a connected component of the boundary $\partial \Sigma$ and $\chi(\Sigma_i)$
is the Euler number of $\Sigma_i$.
This is the combinatorial scalar curvature for the generalized tangential sphere packing metrics introduced in \cite{XZ}.
\end{remark}

We have the following local rigidity for the generalized Thurston's sphere packing metrics.

\begin{theorem}\label{Thm: local rigidity}
Suppose $(\Sigma,\mathcal{T})$ is an ideally triangulated compact 3-dimensional manifold with boundary.
Let $\Phi: E\rightarrow [0,\mathrm{arccos}\frac{1}{3}]$ be a weight defined on the edges of $(\Sigma,\mathcal{T})$.
If $\overline{r}$ is a non-degenerate generalized Thurston's sphere packing metric on $(\Sigma,\mathcal{T},\Phi)$ with $K_{ij}(\overline{r})\sin^2\Phi_{ij}\geq0$ for any edge $\{ij\}\in E$,
then there exists a neighborhood $U$ of $\overline{r}$ such that the generalized Thurston's sphere packing metric in $U$ is determined by its combinatorial scalar curvature.
\end{theorem}

\begin{remark}
As mentioned in \cite{HX}, there exists a subtle difference between the notion of the infinitesimal rigidity in Theorem \ref{Thm: infinitesimal rigidity} and the notion of the local rigidity in Theorem \ref{Thm: local rigidity}.
In fact, the local rigidity implies the infinitesimal rigidity, but the reverse is not true due to dimensional differences.
Note that the condition $K_{ij}\sin^2\Phi_{ij}\geq0$ in Theorem \ref{Thm: local rigidity} includes the case that $K_{ij}(r)<0$ with $\Phi_{ij}=0$ for $\{ij\}\in E$.
Specially, if $\Phi\equiv0$, then Theorem \ref{Thm: local rigidity} is reduced to the local rigidity of the generalized tangential sphere packing.
And the global rigidity of the generalized tangential sphere packings was proved in \cite{XZ}.
\end{remark}

As a direct consequence of Theorem \ref{Thm: local rigidity}, we have the following corollary.
\begin{corollary}
Suppose $(\Sigma,\mathcal{T})$ is an ideally triangulated compact 3-dimensional manifold with boundary.
Let $\Phi: E\rightarrow [0,\mathrm{arccos}\frac{1}{3}]$ be a weight defined on the edges of $(\Sigma,\mathcal{T})$.
If $\overline{r}$ is a non-degenerate generalized Thurston's sphere packing metric on $(\Sigma,\mathcal{T},\Phi)$ with $K_{ij}(\overline{r})\geq0$ for any edge $\{ij\}\in E$,
then there exists a neighborhood $U$ of $\overline{r}$ such that the generalized Thurston's sphere packing metric in $U$ is determined by its combinatorial scalar curvature.
\end{corollary}

The paper is organized as follows.
In Section \ref{section 2}, we prove the simply connectedness of the admissible space of the generalized Thurston's sphere packing metrics for a truncated tetrahedron, i.e., Theorem \ref{Thm: connected}.
In Section \ref{section 3}, we prove Theorem \ref{Thm: infinitesimal rigidity} and Theorem \ref{Thm: local rigidity}.
~\\
\\
\textbf{Acknowledgements}\\[8pt]
The first author thanks Professor Tian Yang at Texas A\&M University for helpful communications.

\section{The admissible space of generalized Thurston's sphere packing metrics on a truncated tetrahedron}\label{section 2}
In this section, we give an explicit characterization of the admissible space of generalized Thurston's sphere packing metrics on a truncated tetrahedron by solving the global version of triangle inequalities for the vertex triangles of a truncated tetrahedron.
The admissible space is proved to be homotopy equivalent to $\mathbb{R}^4_{>0}$ and therefore simply connected.
The approach is taken from He-Xu's work on Thurston's sphere packing metrics for a Euclidean tetrahedron \cite{HX}.
Similar approaches are also used in \cite{Xu1,Xu2,Xu3} and others.

Suppose $\sigma=\{ijkh\}$ is a truncated tetrahedron with four vertex triangles $\triangle_\nu, \nu=i,j,k,h$ as shown in Figure \ref{Figure 2}.
The admissible space $\Omega_{ijkh}$ of the generalized Thurston's sphere packing metrics on $\sigma=\{ijkh\}$ is the set of $(r_i, r_j, r_k, r_h)\in \mathbb{R}^4_{>0}$ such that the hyper-ideal tetrahedron with edge lengths $l_{ij}$,\ $l_{ik}$,\ $l_{ih}$,\ $l_{jk}$,\ $l_{jh}$,\ $l_{kh}$ is non-degenerate.
For a non-degenerate hyper-ideal tetrahedron, Luo-Yang \cite{L-Y} gave the following characterization.
\begin{proposition}(\cite{L-Y}, Proposition 4.4)\label{Prop: L-Y}
The truncated tetrahedron $\sigma=\{ijkh\}$ with edge lengths $(l_{ij}, l_{ik}, l_{ih}, l_{jk}, l_{jh}, l_{kh})\in \mathbb{R}^6_{>0}$ determines a non-degenerate hyper-ideal tetrahedron if and only if the four vertex triangles $\triangle_\nu, \nu=i,j,k,h,$ are non-degenerate, i.e., the triangle inequalities hold for any vertex triangle $\triangle_\nu$.
\end{proposition}

We first consider the vertex triangle $\triangle_i$ with edge lengths $x^i_{jk}, x^i_{jh}, x^i_{kh}$.
To simplify the notations, we set
\begin{equation*}
\begin{aligned}
&\sinh l_{ij}=s_{ij},\ \cosh l_{ij}=c_{ij},\ \tanh r_i=t_i,\ \\
&\cos\Phi_{ij}=a,\ \cos\Phi_{ik}=b,\ \cos\Phi_{ih}=c,\ \\ &\cos\Phi_{jk}=d,\ \cos\Phi_{jh}=e,\  \cos\Phi_{kh}=f.
\end{aligned}
\end{equation*}
Some of the following formulas have been calculated in \cite{XZ}. For completeness, we include them here.
By the cosine law for a right-angled hyperbolic hexagon, we have
\begin{equation}\label{Eq: F11}
\begin{aligned}
\cosh x^i_{jk}=& \frac{\cosh l_{ij}\cosh l_{ik}+\cosh l_{jk}}{\sinh l_{ij}\sinh l_{ik}}=\frac{c_{ij}c_{ik}+c_{jk}}{s_{ij}s_{ik}},\\
\cosh x^i_{jh}=& \frac{\cosh l_{ij}\cosh l_{ih}+\cosh l_{jh}}{\sinh l_{ij}\sinh l_{ih}}=\frac{c_{ij}c_{ih}+c_{jh}}{s_{ij}s_{ih}},\\
\cosh x^i_{kh}=& \frac{\cosh l_{ik}\cosh l_{ih}+\cosh l_{kh}}{\sinh l_{ik}\sinh l_{ih}}=\frac{c_{ik}c_{ih}+c_{kh}}{s_{ik}s_{ih}}.
\end{aligned}
\end{equation}
Note that $x^i_{jk}, x^i_{jh}, x^i_{kh}$ satisfy the triangle inequality if and only if
\begin{equation*}
\sinh \frac{x^i_{jk}+x^i_{jh}+x^i_{kh}}{2}\sinh \frac{x^i_{jk}+x^i_{jh}-x^i_{kh}}{2}\sinh \frac{x^i_{jk}-x^i_{jh}+x^i_{kh}}{2}\sinh \frac{-x^i_{jk}+x^i_{jh}+x^i_{kh}}{2}>0.
\end{equation*}
By direct calculations, we have
\begin{equation*}
\begin{aligned}
&4\sinh \frac{x^i_{jk}+x^i_{jh}+x^i_{kh}}{2}\sinh \frac{x^i_{jk}+x^i_{jh}-x^i_{kh}}{2}\sinh \frac{x^i_{jk}-x^i_{jh}+x^i_{kh}}{2}\sinh \frac{-x^i_{jk}+x^i_{jh}+x^i_{kh}}{2}\\
=&1+2\cosh x^i_{jk}\cosh x^i_{jh}\cosh x^i_{kh}-\cosh^2 x^i_{jk}-\cosh^2x^i_{jh}-\cosh^2x^i_{kh}\\
=&1+2\frac{c_{ij}c_{ik}+c_{jk}}{s_{ij}s_{ik}}\cdot \frac{c_{ij}c_{ih}+c_{jh}}{s_{ij}s_{ih}}\cdot \frac{c_{ik}c_{ih}+c_{kh}}{s_{ik}s_{ih}}
-(\frac{c_{ij}c_{ik}+c_{jk}}{s_{ij}s_{ik}})^2\\
&-(\frac{c_{ij}c_{ih}+c_{jh}}{s_{ij}s_{ih}})^2
-(\frac{c_{ik}c_{ih}+c_{kh}}{s_{ik}s_{ih}})^2\\
=&\frac{1}{s^2_{ij}s^2_{ik}s^2_{ih}}[s^2_{ij}s^2_{ik}s^2_{ih}+2(c_{ij}c_{ik}+c_{jk})(c_{ij}c_{ih}+c_{jh})(c_{ik}c_{ih}+c_{kh})\\
   & \ \ \ \ \ \ \ \ \ \ \ \ \ \ \
   -(c_{ij}c_{ik}+c_{jk})^2s_{ih}^2-(c_{ij}c_{ih}+c_{jh})^2s_{ik}^2-(c_{ik}c_{ih}+c_{kh})^2s_{ij}^2],
\end{aligned}
\end{equation*}
where (\ref{Eq: F11}) is used in the last line.
Furthermore, using the formula $s^2_{ij}=c_{ij}^2-1$, we have
\begin{equation}\label{Eq: Q_1}
\begin{aligned}
Q_1:=&s^2_{ij}s^2_{ik}s^2_{ih}+2(c_{ij}c_{ik}+c_{jk})(c_{ij}c_{ih}+c_{jh})(c_{ik}c_{ih}+c_{kh})\\
   &-(c_{ij}c_{ik}+c_{jk})^2s_{ih}^2-(c_{ij}c_{ih}+c_{jh})^2s_{ik}^2-(c_{ik}c_{ih}+c_{kh})^2s_{ij}^2\\
=&(c_{ij}^2-1)(c_{ik}^2-1)(c_{ih}^2-1)+2(c_{ij}c_{ik}+c_{jk})(c_{ij}c_{ih}+c_{jh})(c_{ik}c_{ih}+c_{kh})\\
 &-(c_{ij}c_{ik}+c_{jk})^2(c_{ih}^2-1)-(c_{ij}c_{ih}+c_{jh})^2(c_{ik}^2-1)-(c_{ik}c_{ih}+c_{kh})^2(c_{ij}^2-1)\\
=&c_{ij}^2+c_{ik}^2+c_{ih}^2+c_{jk}^2+c_{jh}^2+c_{kh}^2
  -c_{ij}^2c_{kh}^2-c_{ik}^2c_{jh}^2-c_{ih}^2c_{jk}^2\\
 &+2c_{ij}c_{ik}c_{jk}+2c_{ik}c_{ih}c_{kh}+2c_{jk}c_{jh}c_{kh}+2c_{ij}c_{ih}c_{jh}\\
 &+2c_{ik}c_{ih}c_{jk}c_{jh}+2c_{ij}c_{ik}c_{jh}c_{kh}+2c_{ij}c_{ih}c_{jk}c_{kh}-1.
\end{aligned}
\end{equation}
Substituting (\ref{Eq: SP metric}) into (\ref{Eq: Q_1}) gives
\begin{equation}\label{Eq: Q_2}
\begin{aligned}
Q_2
=&\frac{1}{\cosh^2 r_i\cosh^2 r_j\cosh^2 r_k\cosh^2 r_h}Q_1\\
=&t_i^2(1-d^2-e^2-f^2-2def)+t_j^2(1-b^2-c^2-f^2-2bcf)\\
&+t_k^2(1-a^2-c^2-e^2-2ace)+t_h^2(1-a^2-b^2-d^2-2abd)\\
&+2t_it_j[(1-f^2)a+bd+ce+bef+cdf]\\
&+2t_it_k[(1-e^2)b+ad+cf+aef+cde]\\
&+2t_it_h[(1-d^2)c+ae+bf+adf+bde]\\
&+2t_jt_k[(1-c^2)d+ab+ef+acf+bce]\\
&+2t_jt_h[(1-b^2)e+ac+df+abf+bcd]\\
&+2t_kt_h[(1-a^2)f+bc+de+abe+acd]\\
&+2abef+2acdf+2bcde+2abd+2ace+2bcf+2def\\
&+a^2+b^2+c^2+d^2+e^2+f^2-a^2f^2-b^2e^2-c^2d^2-1.
\end{aligned}
\end{equation}

Note that $Q_2$ is symmetric in $i, j, k, h$,
 hence the triangle inequalities hold for $\triangle_i, \triangle_j, \triangle_k, \triangle_h$ if and only if $Q_2>0$.
In summary, we have the following result.

\begin{lemma}\label{Lem: non-degenerate}
A hyper-ideal tetrahedron $\sigma=\{ijkh\}$ generated by a generalized Thurston's sphere packing metric is non-degenerate in $\mathbb{H}^3$ if and only if $Q_2>0$.
\end{lemma}

\begin{figure}[!ht]
\centering
\begin{overpic}[scale=0.5]{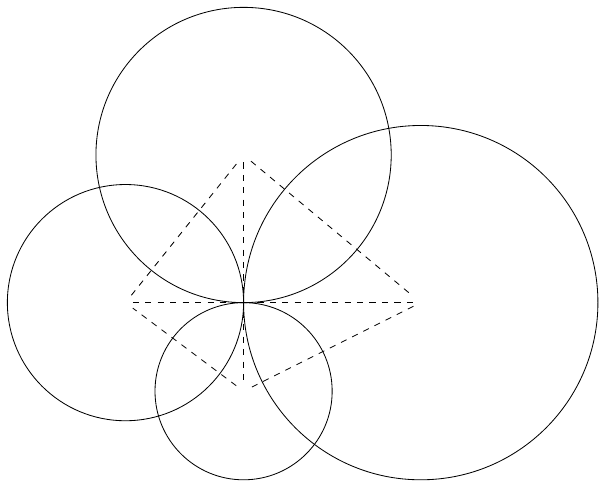}
\put (16,24) {$i$}
\put (37,56) {$h$}
\put (37,8) {$k$}
\put (72,24) {$j$}
\end{overpic}
\caption{Degenerate sphere packing metric}
\label{Figure: four circles}
\end{figure}
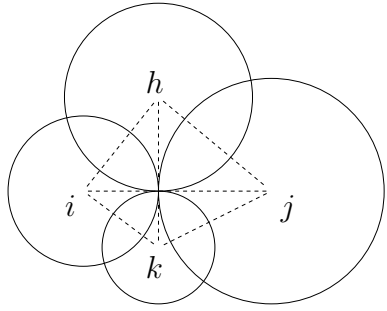

If we assume the weight $\Phi: E\rightarrow [0,\frac{\pi}{2}]$ on $(\Sigma,\mathcal{T})$,
then there are some trivial cases that, for any $r=(r_i, r_j,r_k, r_h)\in \mathbb{R}^4_{>0}$, the hyper-ideal tetrahedron generated by the generalized Thurston's sphere packing metric defined by (\ref{Eq: SP metric}) degenerates.
For example, if the intersection angles on the hyper-ideal tetrahedron satisfies
$a=f=1$ and  $b=c=d=e=0$,
then it is direct to check that $Q_2=0$ for any $r=(r_i, r_j,r_k, r_h)\in \mathbb{R}^4_{>0}$,
which implies that the hyper-ideal tetrahedron degenerates by Lemma \ref{Lem: non-degenerate}.
Please refer to Figure \ref{Figure: four circles} for a configuration of the four spheres.
One can also choose another pair of opposite edges with weight $1$ and the other edges with weight $0$, in
which case the hyper-ideal tetrahedron also degenerates for any $r=(r_i, r_j,r_k, r_h)\in \mathbb{R}^4_{>0}$.

Let $S$ be the set of weights on the edges of the hyper-ideal tetrahedron with value $1$ on a pair of opposite edges and value $0$ on the other edges, i.e.
\begin{eqnarray}\label{Eq: weight S}
\begin{aligned}
S=&\{(\Phi_{ij},\Phi_{ik},\Phi_{ih},\Phi_{jk},\Phi_{jh},\Phi_{kh})\in [0,\frac{\pi}{2}]^6| \\
&(0,\frac{\pi}{2},\frac{\pi}{2},\frac{\pi}{2},\frac{\pi}{2},0), (\frac{\pi}{2},0,\frac{\pi}{2},\frac{\pi}{2},0,\frac{\pi}{2}), (\frac{\pi}{2},\frac{\pi}{2},0,0,\frac{\pi}{2},\frac{\pi}{2})\}.
\end{aligned}
\end{eqnarray}

Set
\begin{equation}\label{Eq: h}
\begin{aligned}
h_i=&t_i(1-d^2-e^2-f^2-2def)+t_j[(1-f^2)a+bd+ce+bef+cdf]\\
&+t_k[(1-e^2)b+ad+cf+aef+cde]+t_h[(1-d^2)c+ae+bf+adf+bde],\\
h_j=&t_j(1-b^2-c^2-f^2-2bcf)+t_i[(1-f^2)a+bd+ce+bef+cdf]\\
&+t_k[(1-c^2)d+ab+ef+acf+bce]+t_h[(1-b^2)e+ac+df+abf+bcd],\\
h_k=&t_k(1-a^2-c^2-e^2-2ace)+t_i[(1-e^2)b+ad+cf+aef+cde]\\
&+t_j[(1-c^2)d+ab+ef+acf+bce]+t_h[(1-a^2)f+bc+de+abe+acd],\\
h_h=&t_h(1-a^2-b^2-d^2-2abd)+t_i[(1-d^2)c+ae+bf+adf+bde]\\
&+t_j[(1-b^2)e+ac+df+abf+bcd]+t_k[(1-a^2)f+bc+de+abe+acd].
\end{aligned}
\end{equation}
Then
\begin{equation}\label{Eq: F5}
Q_2=t_ih_i+t_jh_j+t_kh_k+t_hh_h+Q_3,
\end{equation}
where
\begin{equation*}
\begin{aligned}
Q_3=&2abef+2acdf+2bcde+2abd+2ace+2bcf+2def\\
&+a^2+b^2+c^2+d^2+e^2+f^2-a^2f^2-b^2e^2-c^2d^2-1.
\end{aligned}
\end{equation*}
Set
\begin{equation*}
\begin{aligned}
A_i=d^2+e^2+f^2+2def-1,\ A_j=b^2+c^2+f^2+2bcf-1,\ \\
A_k=a^2+c^2+e^2+2ace-1,\ A_h=a^2+b^2+d^2+2abd-1.
\end{aligned}
\end{equation*}

\begin{lemma}\label{Lem: Q_3>0}
If there exists at least one $\nu\in \{i,j,k,h\}$ such that $A_\nu\geq0$, then $Q_3\geq0$.
In particular, if $A_\nu>0$, then $Q_3>0$.
\end{lemma}
\proof
Note that $Q_3$ is increasing in $a,b,c,d,e,f\in [0,1]$.
If $A_i=d^2+e^2+f^2+2def-1\geq0$, then
\begin{equation*}
\begin{aligned}
Q_3=&(d^2+e^2+f^2+2def-1)+a^2(1-f^2)+b^2(1-e^2)+c^2(1-d^2)\\
&+2abef+2acdf+2bcde+2abd+2ace+2bcf\\
\geq&0.
\end{aligned}
\end{equation*}
Specially, if $A_i>0$, then $Q_3>0$.
Similar results apply to $A_j, A_k, A_h$.
\qed

By Lemma \ref{Lem: non-degenerate},  $r=(r_i,r_j,r_k,r_h)\in \mathbb{R}_{>0}^4$ is a degenerate generalized Thurston's sphere packing metric if and only if $t=(t_i,t_j,t_k,t_h)=(\tanh r_i, \tanh r_j, \tanh r_k, \tanh r_h)\in (0,1)^4$ satisfies
\begin{equation}\label{Eq: F6}
Q_2=t_ih_i+t_jh_j+t_kh_k+t_hh_h+Q_3\leq0,
\end{equation}
which is equivalent to
\begin{equation*}
A_it^2_i+B_it_i+C_i\geq0
\end{equation*}
with
\begin{equation}\label{Eq: F7}
\begin{aligned}
A_i=& d^2+e^2+f^2+2def-1,\\
B_i=&-\big[2t_j\big((1-f^2)a+bd+ce+bef+cdf\big)\\
&\ \ \ \
+2t_k\big((1-e^2)b+ad+cf+aef+cde\big)\\
&\ \ \ \
+2t_h\big((1-d^2)c+ae+bf+adf+bde\big)\big],\\
C_i=&-\big[t_j^2\big(1-b^2-c^2-f^2-2bcf\big)+t_k^2\big(1-a^2-c^2-e^2-2ace\big)\\
&\ \ \ \
+t_h^2\big(1-a^2-b^2-d^2-2abd\big)+2t_jt_k\big((1-c^2)d+ab+ef+acf+bce\big)\\
&\ \ \ \
+2t_jt_h\big((1-b^2)e+ac+df+abf+bcd\big)\\
&\ \ \ \
+2t_kt_h\big((1-a^2)f+bc+de+abe+acd\big)+Q_3\big].
\end{aligned}
\end{equation}

\begin{lemma}\label{Lem: Delta_i}
If $A_i=d^2+e^2+f^2+2def-1>0$, then $\Delta_i=B^2_i-4A_iC_i>0$.
\end{lemma}
\proof
Lengthy and direct calculations give
\begin{equation}\label{Eq: F1}
\begin{aligned}
\Delta_i=&-4[t_j^2(1-f^2)+t_k^2(1-e^2)+t_h^2(1-d^2)\\
&\ \ \ \ \
+2t_jt_k(d+ef)+2t_jt_h(e+df)+2t_kt_h(f+de)]\\
&\times(1-a^2-b^2-c^2-d^2-e^2-f^2+a^2f^2+b^2e^2+c^2d^2\\
&\ \ \ \ \
-2ace-2abd-2bcf-2def-2abef-2acdf-2bcde)\\
&+4(d^2+e^2+f^2+2def-1)Q_3\\
=&4Q_3[(d^2+e^2+f^2+2def-1)+t_j^2(1-f^2)+t_k^2(1-e^2)+t_h^2(1-d^2)\\
&+2t_jt_k(d+ef)+2t_jt_h(e+df)+2t_kt_h(f+de)].
\end{aligned}
\end{equation}
If $A_i>0$, then $Q_3>0$ by Lemma \ref{Lem: Q_3>0}.
This implies $\Delta_i>0$ by (\ref{Eq: F1}).
\qed

\begin{remark}
If we take $Q_2$ as a quadratic function of $t_\nu$ for $\nu\in\{j,k,h\}$, results similar to Lemma \ref{Lem: Delta_i} hold for $\nu\in\{j,k,h\}$.
\end{remark}


\begin{lemma}\label{Lem: two non-positive}
Let $\sigma=\{ijkh\}$ be a truncated tetrahedron with the weight $\Phi\in [0, \frac{\pi}{2}]^6\setminus S$.
If $r=(r_i,r_j,r_k,r_h)\in \mathbb{R}_{>0}^4$ is a degenerate generalized Thurston's sphere packing metric,
then there exists no subset $\{\mu,\nu\}\subseteq\{i,j,k,h\}$ such that $h_\mu\leq0$ and $h_\nu\leq0$.
\end{lemma}
\proof
Without loss of generality, we assume $h_i\leq0$ and $h_j\leq0$.
By (\ref{Eq: h}), we have
\begin{equation}\label{Eq: F8}
\begin{aligned}
t_i(d^2+e^2+f^2+2def-1)\geq &t_j[(1-f^2)a+bd+ce+bef+cdf]\\
&+t_k[(1-e^2)b+ad+cf+aef+cde]\\
&+t_h[(1-d^2)c+ae+bf+adf+bde]
\end{aligned}
\end{equation}
and
\begin{equation}\label{Eq: F9}
\begin{aligned}
t_j(b^2+c^2+f^2+2bcf-1)\geq &t_i[(1-f^2)a+bd+ce+bef+cdf]\\
&+t_k[(1-c^2)d+ab+ef+acf+bce]\\
&+t_h[(1-b^2)e+ac+df+abf+bcd].
\end{aligned}
\end{equation}
This implies
\begin{equation*}
A_i=d^2+e^2+f^2+2def-1\geq0, \quad A_j=b^2+c^2+f^2+2bcf-1\geq0,
\end{equation*}
and
\begin{equation*}
\begin{aligned}
M:=&(d^2+e^2+f^2+2def-1)(b^2+c^2+f^2+2bcf-1)\\
&-[(1-f^2)a+bd+ce+bef+cdf]^2\\
\geq&0.
\end{aligned}
\end{equation*}
However, direct calculations give
\begin{equation*}
M=-(1-f^2)Q_3\leq0
\end{equation*}
by Lemma \ref{Lem: Q_3>0}.
Therefore, $M\equiv0$, which implies $f=1$ or $Q_3=0$.

In the case of $f=1$,
since $M\equiv0$, then the inequality (\ref{Eq: F8}) should be an equality and the second and third terms in the righthand side of (\ref{Eq: F8}) should be $0$.
This  implies $b=c=0$ and $a(d+e)=0$.
Similarly, we have $d=e=0$ and $a(b+c)=0$ by (\ref{Eq: F9}).
By the condition that the weight $\Phi\in [0, \frac{\pi}{2}]^6\setminus S$, we have $a\neq 1$.
Substituting $b=c=d=e=0, f=1, a\neq1$ into (\ref{Eq: Q_2}) gives
\begin{equation*}
Q_2=t_k^2(1-a^2)+t_h^2(1-a^2)+2t_kt_h(1-a^2)>0.
\end{equation*}
This contradicts with the condition that $r$ is a degenerate generalized Thurston's sphere packing metric.

In the case of $Q_3=0$, we have $A_\nu\leq0$ for any $\nu\in\{i,j,k,h\}$ by Lemma \ref{Lem: Q_3>0}.
By $Q_2\leq0$ and (\ref{Eq: Q_2}), we have $A_\nu\equiv0$ for any $\nu\in\{i,j,k,h\}$ ,
which implies $Q_3>0$ due to the weight $\Phi\in [0, \frac{\pi}{2}]^6\setminus S$.
This is a contradiction.
\qed

By Lemma \ref{Lem: two non-positive}, if $r=(r_i,r_j,r_k,r_h)\in \mathbb{R}_{>0}^4$ is a degenerate generalized Thurston's sphere packing metric on a hyper-ideal tetrahedron with the weight $\Phi\in [0, \frac{\pi}{2}]^6\setminus S$, then there are only two cases.
One is that $h_\nu>0$ for all $\nu\in \{i,j,k,h\}$, the other is that one of $h_i, h_j, h_k, h_h$ is negative and the others are positive.
To remove the former, we need to add some restrictions on the weight $\Phi$.

\begin{lemma}\label{Lem: four positive}
Let $\sigma=\{ijkh\}$ be a truncated tetrahedron with the weight $\Phi: E\rightarrow [0,\mathrm{arccos}\frac{1}{3}]$.
If $r=(r_i,r_j,r_k,r_h)\in \mathbb{R}_{>0}^4$ is a degenerate generalized Thurston's sphere packing metric,
then it is impossible that $h_i>0,\ h_j>0,\ h_k>0$ and $h_h>0$.
\end{lemma}
\proof
Assume $h_i>0,\ h_j>0,\ h_k>0$ and $h_h>0$.
Since $h_i>0$, then
\begin{equation*}
\begin{aligned}
t_i&(d^2+e^2+f^2+2def-1)<t_j[(1-f^2)a+bd+ce+bef+cdf]\\
&+t_k[(1-e^2)b+ad+cf+aef+cde]+t_h[(1-d^2)c+ae+bf+adf+bde].
\end{aligned}
\end{equation*}
If $A_i=d^2+e^2+f^2+2def-1\geq0$, then $Q_3\geq0$ by Lemma \ref{Lem: Q_3>0}.
This implies $Q_2>0$ by (\ref{Eq: F5}),
which contradicts with the condition that $r$ is a degenerate generalized Thurston's sphere packing metric by Lemma \ref{Lem: non-degenerate}.
Therefore, we only need to consider the case that $A_i<0$,\ $A_j<0$,\ $A_k<0$ and $A_h<0$.
Note that $\Phi: E\rightarrow [0,\mathrm{arccos}\frac{1}{3}]$, i.e., $a,b,c,d,e,f\in [\frac{1}{3},1]$.
Since $Q_3$ is increasing in $a,b,c,d,e,f$, then $Q_3\geq0$ and $Q_3=0$ if and only if $a=b=c=d=e=f=\frac{1}{3}$.
This implies $Q_2>0$ by (\ref{Eq: F5}),
which contradicts with the condition that $r$ is a degenerate generalized Thurston's sphere packing metric by Lemma \ref{Lem: non-degenerate}.
\qed

\begin{remark}\label{Rmk: key remark}
By Lemma \ref{Lem: Q_3>0}, if $A_\nu\geq0$ for some $\nu\in \{i,j,k,h\}$, then $Q_3\geq0$.
However, if $A_\nu<0$ for all $\nu\in \{i,j,k,h\}$, then the sign of $Q_3$ is unclear, neither is $Q_2$ by (\ref{Eq: F5}).
To ensure $Q_3\geq0$, we add a strong condition that $\Phi: E\rightarrow [0,\mathrm{arccos}\frac{1}{3}]$ on the weight.
In this case, the set $S$ defined by (\ref{Eq: weight S}) is empty.
And Lemma \ref{Lem: Q_3>0}, Lemma \ref{Lem: Delta_i} and Lemma \ref{Lem: two non-positive} still hold.
\end{remark}

Combining Lemma \ref{Lem: two non-positive} with Lemma \ref{Lem: four positive}, we have the following corollary.

\begin{corollary}\label{Cor: one negative and three positive}
Let $\sigma=\{ijkh\}$ be a truncated tetrahedron with the weight $\Phi: E\rightarrow [0,\mathrm{arccos}\frac{1}{3}]$.
If $r=(r_i,r_j,r_k,r_h)\in \mathbb{R}_{>0}^4$ is a degenerate generalized Thurston's sphere packing metric,
then one of $h_i, h_j, h_k, h_h$ is negative and the others are positive.
\end{corollary}

Now we can prove the main result of this section.
\begin{theorem}\label{Thm: connected}
Let $\sigma=\{ijkh\}$ be a truncated tetrahedron with the weight $\Phi: E\rightarrow [0,\mathrm{arccos}\frac{1}{3}]$.
Then the admissible space $\Omega_{ijkh}(\Phi)$ of the generalized Thurston's sphere packing metric $r=(r_i,r_j,r_k,r_h)\in \mathbb{R}_{>0}^4$ is a simply connected non-empty open set.
\end{theorem}
\proof
Suppose $r=(r_i,r_j,r_k,r_h)\in \mathbb{R}_{>0}^4$ is a degenerate generalized Thurston's sphere packing metric, then $t=(t_i,t_j,t_k,t_h)=(\tanh r_i, \tanh r_j, \tanh r_k, \tanh r_h)\in (0,1)^4$ satisfies (\ref{Eq: F6}).
By Corollary \ref{Cor: one negative and three positive}, one of $h_i, h_j, h_k, h_h$ is negative and the others are positive.
Without loss of generality, we assume $h_i<0$.
Then $A_i=d^2+e^2+f^2+2def-1>0$ by the definition of $h_i$.
Take $Q_2\leq0$ as a quadratic inequality $A_it^2_i+B_it_i+C_i\geq0$, where $A_i, B_i, C_i$ are
given by (\ref{Eq: F7}).
As $A_i>0$, we have $\Delta_i=B^2_i-4A_iC_i>0$ by Lemma \ref{Lem: Delta_i}, which implies
\begin{equation*}
t_i\geq\frac{-B_i+\sqrt{\Delta_i}}{2A_i} \quad \text{or} \quad t_i\leq\frac{-B_i-\sqrt{\Delta_i}}{2A_i}.
\end{equation*}
The definition of $h_i$ in (\ref{Eq: h}) implies $h_i=-\frac{1}{2}(2A_it_i+B_i)$.
Then $h_i<0$ is equivalent to $t_i>-\frac{B_i}{2A_i}$,
which implies $t_i\geq\frac{-B_i+\sqrt{\Delta_i}}{2A_i}$.

In the case that $A_i>0$, we set
\begin{equation}\label{Eq: V_i}
V_i=\{r=(r_i,r_j,r_k,r_h)\in \mathbb{R}_{>0}^4|t_i\geq\frac{-B_i+\sqrt{\Delta_i}}{2A_i}\},
\end{equation}
which is bounded by an analytical function defined on $\mathbb{R}_{>0}^3$.
For any $(r_i,r_j,r_k,r_h)\in V_i$, it is direct to check that $A_it^2_i+B_it_i+C_i\geq0$, which implies $Q_2\leq 0$, thus $V_i\subseteq \mathbb{R}^4\setminus \Omega_{ijkh}(\Phi)$.
One can define $V_j,V_k,V_h$ similarly, if the corresponding $A_j>0, A_k>0, A_h>0$.
Similar arguments imply $V_j, V_k, V_h\subseteq \mathbb{R}^4\setminus \Omega_{ijkh}(\Phi)$.

If $A_\mu>0$ for $\mu\in P\subseteq \{i,j,k,h\}$ and $A_\omega\leq0$ for $\omega\in \{i,j,k,h\}\backslash P$,
then the space of degenerate generalized Thurston's sphere packing metrics is $\bigcup_{\mu\in P}V_\mu$.
Therefore,
\begin{equation*}
\Omega_{ijkh}(\Phi)
=\mathbb{R}_{>0}^4\setminus\bigcup_{\mu\in P}{V_\mu}.
\end{equation*}

For any $r\in V_i$, we have $t_i>-\frac{B_i}{2A_i}$, which is equivalent to $h_i<0$.
Similarly, for $r\in V_j$, we have $h_j<0$.
Then Corollary \ref{Cor: one negative and three positive} implies $V_i\bigcap V_j=\varnothing$.
Similarly, we have $V_i,V_j,V_k,V_h$ are mutually disjoint, if they are non-empty.
Therefore, $\Omega_{ijkh}(\Phi)$ is homotopy equivalent to $\mathbb{R}_{>0}^4$ and hence is a simply connected non-empty open set.
\qed

By the proof of Theorem \ref{Thm: connected}, the admissible space $\Omega_{ijkh}(\Phi)=\mathbb{R}_{>0}^4$
if and only if $A_i\leq0\ ,A_j\leq0,\ A_k\leq0,\ A_h\leq0$.
Furthermore, we have the following results.

\begin{corollary}\label{Cor: admissible space 1}
Let $\sigma=\{ijkh\}$ be a truncated tetrahedron with the weight $\Phi: E\rightarrow [0,\mathrm{arccos}\frac{1}{3}]$.
If any three of $A_i,A_j,A_k,A_h$ are non-positive,
then the admissible space $\Omega_{ijkh}(\Phi)=\mathbb{R}_{>0}^4$.
Specially, if any three of the following four formulas hold,
\begin{align*}
&\Phi_{ij}+\Phi_{jk}+\Phi_{ik}\geq\pi,\
\Phi_{ik}+\Phi_{ih}+\Phi_{kh}\geq\pi,\ \\
&\Phi_{ij}+\Phi_{ih}+\Phi_{jh}\geq\pi,\
\Phi_{ij}+\Phi_{ik}+\Phi_{jk}\geq\pi,
\end{align*}
then the admissible space $\Omega_{ijkh}(\Phi)=\mathbb{R}_{>0}^4$.
\end{corollary}
\proof
Without loss of generality, we assume $A_j\leq0, A_k\leq0, A_h\leq0$.
By the definition of $h_j,h_k,h_h$ in (\ref{Eq: h}) and Lemma \ref{Lem: two non-positive},
we have $h_j>0,\ h_k>0,\ h_h>0$.
If $r=(r_i,r_j,r_k,r_h)\in \mathbb{R}_{>0}^4$ is a degenerate generalized Thurston's sphere packing metric,
then $h_i<0$ by Corollary \ref{Cor: one negative and three positive}.
Then $A_i=d^2+e^2+f^2+2def-1>0$ and hence $Q_3>0$ by Lemma \ref{Lem: Q_3>0}.
Note that $\Phi\in [0,\mathrm{arccos}\frac{1}{3}]^6$, i.e., $a,b,c,d,e,f\in [\frac{1}{3},1]$, then $Q_3-A_i>0$. This implies
\begin{equation*}
-B_i+\sqrt{\Delta_i}-2A_i
>2\sqrt{A_iQ_3}-2A_i
=2\sqrt{A_i}\cdot\frac{Q_3-A_i}{\sqrt{Q_3}+\sqrt{A_i}}>0.
\end{equation*}
Hence, $\frac{-B_i+\sqrt{\Delta_i}}{2A_i}>1$.
Combining with $t_i=\tanh r_i\in (0,1)$ gives $V_i=\emptyset$ by (\ref{Eq: V_i}).
Thus $r=(r_i,r_j,r_k,r_h)\in \mathbb{R}_{>0}^4$ is a non-degenerate generalized Thurston's sphere packing metric and hence $\Omega_{ijkh}(\Phi)=\mathbb{R}_{>0}^4$.

For the second part, note that
\begin{equation*}
\begin{aligned}
A_i=&\cos^2\Phi_{jk}+\cos^2\Phi_{jh}+\cos^2\Phi_{kh}
+2\cos\Phi_{jk}\cos\Phi_{jh}\cos\Phi_{kh}-1\\
=&4\cos\frac{\Phi_{jk}+\Phi_{jh}+\Phi_{kh}}{2}
\cos\frac{\Phi_{jk}+\Phi_{jh}-\Phi_{kh}}{2}
\cos\frac{\Phi_{jk}-\Phi_{jh}+\Phi_{kh}}{2}
\cos\frac{-\Phi_{jk}+\Phi_{jh}+\Phi_{kh}}{2}.
\end{aligned}
\end{equation*}
Since $\Phi_{ij}\in [0, \mathrm{acrcos}\frac{1}{3}]\subseteq[0,\frac{\pi}{2}]$,
then
\begin{equation*}
\cos\frac{\Phi_{jk}+\Phi_{jh}-\Phi_{kh}}{2}\geq0,\
\cos\frac{\Phi_{jk}-\Phi_{jh}+\Phi_{kh}}{2}\geq0,\
\cos\frac{-\Phi_{jk}+\Phi_{jh}+\Phi_{kh}}{2}\geq0.
\end{equation*}
Therefore, if $\Phi_{jk}+\Phi_{jh}+\Phi_{kh}\geq\pi$, then $\cos\frac{\Phi_{jk}+\Phi_{jh}+\Phi_{kh}}{2}\leq0$ and hence $A_i\leq0$.
Similarly, if
\begin{equation*}
\Phi_{ik}+\Phi_{ih}+\Phi_{kh}\geq\pi,\
\Phi_{ij}+\Phi_{ih}+\Phi_{jh}\geq\pi,\
\Phi_{ij}+\Phi_{ik}+\Phi_{jk}\geq\pi,
\end{equation*}
then $A_j\leq0,\ A_k\leq0,\ A_h\leq0$ respectively.
\qed

\begin{corollary}\label{Cor: admissible space 2}
Let $\sigma=\{ijkh\}$ be a truncated tetrahedron with the weight $\Phi: E\rightarrow [0,\mathrm{arccos}\frac{1}{3}]$.
If $\cos\Phi\equiv C\in[\frac{1}{3},1]$,
then the admissible space $\Omega_{ijkh}(\Phi)=\mathbb{R}_{>0}^4$.
\end{corollary}
\proof
If $\cos\Phi\equiv C\in[\frac{1}{3},\frac{1}{2}]$,
then $A_i=d^2+e^2+f^2+2def-1=3C^2+2C^3-1\leq0$.
Similarly, $A_j\leq0,\ A_k\leq0,\ A_h\leq0$.
Then $\Omega_{ijkh}(\Phi)=\mathbb{R}_{>0}^4$ by Corollary \ref{Cor: admissible space 1}.
If $\cos\Phi\equiv C\in[\frac{1}{2},1]$,
then $1-3C^2-2C^3\leq0$.
Note that $t_i=\tanh r_i\in (0,1)$,
then by (\ref{Eq: Q_2}), we have
\begin{equation*}
\begin{aligned}
Q_2
=&(t_i^2+t_j^2+t_k^2+t_h^2)(1-3C^2-2C^3)\\
&+2(t_it_j+t_it_k+t_it_h+t_jt_k+t_jt_h+t_kt_h)(C^3+2C^2+C)\\
&+3C^4+8C^3+6C^2-1\\
\geq&2(t_it_j+t_it_k+t_it_h+t_jt_k+t_jt_h+t_kt_h)(C^3+2C^2+C)\\
&+3C^4-6C^2+3\\
>&0,
\end{aligned}
\end{equation*}
which implies $\Omega_{ijkh}(\Phi)=\mathbb{R}_{>0}^4$ by Lemma \ref{Lem: non-degenerate}.
\qed

\begin{remark}
As a special case that $\Phi\equiv0$, Corollary \ref{Cor: admissible space 2} shows that the admissible space $\Omega_{ijkh}(\Phi)$ for the generalized tangential sphere packing metrics is $\mathbb{R}_{>0}^4$, which was proved in \cite{XZ}.
\end{remark}

\section{Rigidity of generalized Thurston's sphere packings}\label{section 3}
Theorem \ref{Thm: connected} shows that the admissible space $\Omega_{ijkh}(\Phi)$ is simply connected for a weighted truncated tetrahedron $\sigma=\{ijkh\}$.
Hence, the admissible space on  $(\Sigma,\mathcal{T},\Phi)$ can be defined to be $\Omega(\Phi)=\bigcap_{\sigma=\{ijkh\}}\Omega_{ijkh}(\Phi)$. Note that $\Omega(\Phi)$ may be not connected.

Denote the volume of a hyper-ideal tetrahedra $\sigma=\{ijkh\}$ by $vol$, which is considered as a function in the dihedral angles.
By the Schl\"{a}fli formula, we have
\begin{equation}\label{Schlafli formula}
2d vol+\sum_{j\sim i}l_{ij}d\beta_{ij,kh}=0.
\end{equation}
One can refer to \cite{Bonahon,L-Y} for more details on the Schl\"{a}fli formula.
Using the volume function $vol$ on a singer hyper-ideal tetrahedra $\sigma=\{ijkh\}$ , we can define the following function
\begin{equation*}
G=\sum_{j\sim i}K_{ij}l_{ij}-2\sum_{\sigma\in T} vol
\end{equation*}
on $(\Sigma,\mathcal{T},\Phi)$, where $K_{ij}$ is the combinatorial Ricci curvature along the edge $\{ij\}\in E$ defined by (\ref{Eq: CRC}).
Then  by (\ref{Schlafli formula}), we have $dG=\sum_{j\sim i}K_{ij}dl_{ij}$, which implies
\begin{equation}\label{Eq: F10}
\frac{\partial G}{\partial r_i}
=\sum_{j;j\sim i}K_{ij}\frac{\partial l_{ij}}{\partial r_i}
=\sum_{j;j\sim i}K_{ij}\cosh \theta_{ij}
=K_i.
\end{equation}
Furthermore,
\begin{equation*}
\begin{aligned}
\frac{\partial^2 G}{\partial r_i\partial r_j}
=&\sum_{j;j\sim i}\frac{\partial K_{ij}}{\partial r_j}\frac{\partial l_{ij}}{\partial r_i}
+\sum_{j;j\sim i}K_{ij}\frac{\partial^2 l_{ij}}{\partial r_i\partial r_j}\\
=&\sum_{s\sim t}\sum_{u\sim v}\frac{\partial l_{uv}}{\partial r_j}\frac{\partial K_{st}}{\partial l_{uv}}\frac{\partial l_{st}}{\partial r_i}
+\sum_{s\sim t}K_{st}\frac{\partial^2 l_{st}}{\partial r_i\partial r_j}\\
=&-\sum_{\{ijkh\}\in T}\sum_{s\sim t}\sum_{u\sim v}\frac{\partial l_{uv}}{\partial r_j}\frac{\partial \beta_{st}}{\partial l_{uv}}\frac{\partial l_{st}}{\partial r_i}
+\sum_{s\sim t}K_{st}\frac{\partial^2 l_{st}}{\partial r_i\partial r_j},
\end{aligned}
\end{equation*}
where $\beta_{ij}:=\beta_{ij,kh}$ for simplification.
This implies
\begin{equation}\label{Eq: Hess G}
\mathrm{Hess}_r G=-\sum_{\{ijkh\}\in T}\bigg(\frac{\partial l}{\partial r}\bigg)_\sigma\bigg(\frac{\partial \beta}{\partial l}\bigg)_\sigma\bigg(\frac{\partial l}{\partial r}\bigg)_\sigma^T+\sum_{j\sim i}K_{ij}\mathrm{Hess}_r(l_{ij}),
\end{equation}
where
\begin{equation*}
\bigg(\frac{\partial l}{\partial r}\bigg)_\sigma
=\left(
   \begin{array}{cccccc}
    \cosh \theta_{ij} & \cosh \theta_{ik} & \cosh \theta_{ih} & 0 & 0 & 0 \\
     \cosh \theta_{ji} & 0 & 0 & \cosh \theta_{jk} & \cosh \theta_{jh} & 0 \\
     0 & \cosh \theta_{ki} & 0 & \cosh \theta_{kj} & 0 & \cosh \theta_{kh} \\
     0 & 0 & \cosh \theta_{hi} & 0 & \cosh \theta_{hj} & \cosh \theta_{hk}
    \end{array}
 \right)
\end{equation*}
and
\begin{equation*}
\bigg(\frac{\partial \beta}{\partial l}\bigg)_\sigma
=\left(
\begin{array}{cccccc}
\frac{\partial\beta_{ij}}{\partial l_{ij}} & \frac{\partial\beta_{ij}}{\partial l_{ik}} & \frac{\partial\beta_{ij}}{\partial l_{ih}} & \frac{\partial\beta_{ij}}{\partial l_{jk}} & \frac{\partial\beta_{ij}}{\partial l_{jh}} & \frac{\partial\beta_{ij}}{\partial l_{kh}} \\
\frac{\partial\beta_{ik}}{\partial l_{ij}} & \frac{\partial\beta_{ik}}{\partial l_{ik}} & \frac{\partial\beta_{ik}}{\partial l_{ih}} & \frac{\partial\beta_{ik}}{\partial l_{jk}} & \frac{\partial\beta_{ik}}{\partial l_{jh}} & \frac{\partial\beta_{ik}}{\partial l_{kh}} \\
\frac{\partial\beta_{ih}}{\partial l_{ij}} & \frac{\partial\beta_{ih}}{\partial l_{ik}} & \frac{\partial\beta_{ih}}{\partial l_{ih}} & \frac{\partial\beta_{ih}}{\partial l_{jk}} & \frac{\partial\beta_{ih}}{\partial l_{jh}} & \frac{\partial\beta_{ih}}{\partial l_{kh}} \\
\frac{\partial\beta_{jk}}{\partial l_{ij}} & \frac{\partial\beta_{jk}}{\partial l_{ik}} & \frac{\partial\beta_{jk}}{\partial l_{ih}} & \frac{\partial\beta_{jk}}{\partial l_{jk}} & \frac{\partial\beta_{jk}}{\partial l_{jh}} & \frac{\partial\beta_{jk}}{\partial l_{kh}} \\
\frac{\partial\beta_{jh}}{\partial l_{ij}} & \frac{\partial\beta_{jh}}{\partial l_{ik}} & \frac{\partial\beta_{jh}}{\partial l_{ih}} & \frac{\partial\beta_{jh}}{\partial l_{jk}} & \frac{\partial\beta_{jh}}{\partial l_{jh}} & \frac{\partial\beta_{jh}}{\partial l_{kh}} \\
\frac{\partial\beta_{kh}}{\partial l_{ij}} & \frac{\partial\beta_{kh}}{\partial l_{ik}} & \frac{\partial\beta_{kh}}{\partial l_{ih}} & \frac{\partial\beta_{kh}}{\partial l_{jk}} & \frac{\partial\beta_{kh}}{\partial l_{jh}} & \frac{\partial\beta_{kh}}{\partial l_{kh}} \\
\end{array}
\right).
\end{equation*}

It is easy to check that the rank of the metric $(\frac{\partial l}{\partial r})_\sigma$ is 4.
We have the following result on the matrix $(\frac{\partial \beta}{\partial l})_\sigma$.

\begin{theorem}[\cite{Luo}]\label{Thm: positive definite}
For a hyper-ideal tetrahedron $\sigma=\{ijkh\}$ with dihedral angle $\beta_{ij}$ and length $l_{ij}$ at the edge $\{ij\}\in E$, the matrix $(\frac{\partial \beta}{\partial l})$ is symmetric and strictly positive definite on $\mathcal{L}$, where $\mathcal{L}$ is the set of vectors $(l_{ij}, l_{ik}, l_{ih}, l_{jk}, l_{jh}, l_{kh})\in \mathbb{R}^6_{>0}$ such that there exists a non-degenerate hyper-ideal tetrahedron $\sigma=\{ijkh\}$ with $l_{ij}$ as the length of the edge $\{ij\}$.
\end{theorem}

\begin{remark}
The simply connectedness of $\mathcal{L}$ has been proved in \cite{Bao 2,F,FP,Luo}.
\end{remark}

By Theorem \ref{Thm: positive definite}, the first term in the right of (\ref{Eq: Hess G}) is a symmetric and strictly negative definite matrix on the admissible space $\Omega(\Phi)$.
Now we consider the second term in the righthand side of (\ref{Eq: Hess G}).
Note that in the special case of generalized sphere packings, i.e., $\Phi\equiv 0$, the second term equals 0 because of $l_{ij}=r_i+r_j$.
Then $\mathrm{Hess}_r G<0$ on $(\Sigma,\mathcal{T},\Phi)$,
which was proved in \cite{XZ} by direct calculations.
Please refer to \cite{XZ} for more details.
In the general case of $\Phi: E\rightarrow [0, \frac{\pi}{2}]$, we have
\begin{equation*}
\begin{aligned}
\frac{\partial^2 l_{ij}}{\partial r^2_i}
=&\frac{\cosh l_{ij}}{\sinh^3 l_{ij}}[\sinh^2 l_{ij}-(\cosh r_i\sinh r_j+\cos \Phi_{ij}\sinh r_i\cosh r_j)^2]\\
=&\frac{\cosh l_{ij}}{\sinh^3 l_{ij}}[(\sinh r_i\sinh r_j+\cos \Phi_{ij}\cosh r_i\cosh r_j)^2-1\\
&-(\cosh r_i\sinh r_j+\cos \Phi_{ij}\sinh r_i\cosh r_j)^2]\\
=&-\frac{\sin^2\Phi_{ij}}{\sinh^3 l_{ij}}\cosh l_{ij}\cosh^2r_j,
\end{aligned}
\end{equation*}
and
\begin{equation*}
\begin{aligned}
\frac{\partial^2 l_{ij}}{\partial r_i\partial r_j}
=&\frac{1}{\sinh^3 l_{ij}}[(\cosh r_i\cosh r_j+\cos \Phi_{ij}\sinh r_i\sinh r_j)(\cosh^2 l_{ij}-1)\\
&-\cosh l_{ij}(\sinh r_i\cosh r_j+\cos \Phi_{ij}\cosh r_i\sinh r_j)\\
&\times(\cosh r_i\sinh r_j+\cos \Phi_{ij}\sinh r_i\cosh r_j)]\\
=&-\frac{\sin^2\Phi_{ij}}{\sinh^3 l_{ij}}\cosh r_i\cosh r_j.
\end{aligned}
\end{equation*}
Then
\begin{equation*}
\mathrm{Hess}_r(l_{ij})
=-\frac{\sin^2\Phi_{ij}}{\sinh^3 l_{ij}}
\left(
   \begin{array}{cc}
    \cosh l_{ij}\cosh^2r_j & \cosh r_i\cosh r_j  \\
     \cosh r_i\cosh r_j & \cosh l_{ij}\cosh^2r_i  \\
    \end{array}
 \right),
\end{equation*}
which is a symmetric and negative definite matrix if $\Phi_{ij}\not\equiv0$ for any edge $\{ij\}\in E$.
Therefore,
\begin{equation*}
\begin{aligned}
\mathrm{Hess}_r G
=&-\sum_{\{ijkh\}\in T}\bigg(\frac{\partial l}{\partial r}\bigg)_\sigma\bigg(\frac{\partial \beta}{\partial l}\bigg)_\sigma\bigg(\frac{\partial l}{\partial r}\bigg)_\sigma^T\\
&-\sum_{j\sim i}\frac{K_{ij}\sin^2\Phi_{ij}}{\sinh^3 l_{ij}}
\left(
   \begin{array}{cc}
    \cosh l_{ij}\cosh^2r_j & \cosh r_i\cosh r_j  \\
     \cosh r_i\cosh r_j & \cosh l_{ij}\cosh^2r_i  \\
    \end{array}
 \right).
\end{aligned}
\end{equation*}

\noindent\textbf{Proof\ of\ Theorem\ \ref{Thm: local rigidity}}
If $\overline{r}$ is a non-degenerate generalized Thurston's sphere packing metric on $(\Sigma,\mathcal{T},\Phi)$ with $K_{ij}(\overline{r})\sin^2\Phi_{ij}\geq0$ for any edge $\{ij\}\in E$,
then $\mathrm{Hess}_r G(\overline{r})<0$ and hence the function $G$ is a strictly concave function on a convex neighborhood $U$ of $\overline{r}$.

Suppose there exist two different generalized Thurston's sphere packing metrics $r_1$ and $r_2$ in $U$ such that $K_i(r_1)=K_i(r_2)$ for all $i\in V$.
Set
\begin{equation*}
f(t)=G((1-t)r_1+tr_2),\ t\in [0,1].
\end{equation*}
Then $f(t)$ is a $C^1$ smooth concave function of $t\in [0,1]$ with
\begin{equation*}
f'(t)=\sum_{i=1}^{N}\nabla_{r_i} G|_{(1-t)r_1+tr_2}\cdot (r_{2,i}-r_{1,i})
=\sum_{i=1}^{N}K_i|_{(1-t)r_1+tr_2}\cdot (r_{2,i}-r_{1,i}).
\end{equation*}
By the assumption that $K_i(r_1)=K_i(r_2)$, we have $f'(0)=f'(1)=0$,
which implies $f'(t)\equiv0$ by the concavity of $f(t)$.
Then $f''(t)\equiv0$ for $t\in [0, 1]$.
Combining with
\begin{equation*}
f''(t)=(r_2-r_1)^T\cdot \mathrm{Hess}_r G|_{(1-t)r_1+tr_2} \cdot (r_2-r_1)
\end{equation*}
for $t\in [0,1]$ and $\mathrm{Hess}_r G$ is negative definite,
we have $r_1=r_2$.
\qed

For the combinatorial Ricci curvature, we rewrite Theorem \ref{Thm: infinitesimal rigidity} as follows.
\begin{theorem}\label{Thm: infinitesimal rigidity 2}
Suppose $(\Sigma,\mathcal{T})$ is an ideally triangulated compact 3-dimensional manifold with boundary.
Let $\Phi: E\rightarrow [0,\mathrm{arccos}\frac{1}{3}]$ be a weight defined on the edges of $(\Sigma,\mathcal{T})$.
Suppose there exists a neighborhood $U\subseteq\Omega(\Phi)$ of $\overline{r}$ such that if $r\in U$ has the same combinatorial Ricci curvature as $\overline{r}$, then $r=\overline{r}$.
\end{theorem}
\proof
Define the following functional
\begin{equation*}
\widetilde{G}(r_1,r_2,...,r_N)
=G-\sum_{j;j\sim i}K_{ij}(\overline{r})l_{ij}.
\end{equation*}
Then by (\ref{Eq: F10}), we have
\begin{equation*}
\frac{\partial \widetilde{G}}{\partial r_i}
=K_i-\sum_{j;j\sim i}K_{ij}(\overline{r})\frac{\partial l_{ij}}{\partial r_i}
=\sum_{j;j\sim i}(K_{ij}-K_{ij}(\overline{r}))\frac{\partial l_{ij}}{\partial r_i},
\end{equation*}
which implies $\nabla \widetilde{G}(r)=\nabla \widetilde{G}(\overline{r})$ by the assumption that $K_{ij}(r)=K_{ij}(\overline{r})$.
Furthermore,
\begin{equation*}
\begin{aligned}
\mathrm{Hess}_r \widetilde{G}
=&-\sum_{\{ijkh\}\in T}\bigg(\frac{\partial l}{\partial r}\bigg)_\sigma\bigg(\frac{\partial \beta}{\partial l}\bigg)_\sigma\bigg(\frac{\partial l}{\partial r}\bigg)_\sigma^T\\
&-\sum_{j\sim i}(K_{ij}-K_{ij}(\overline{r}))\frac{\sin^2\Phi_{ij}}{\sinh^3 l_{ij}}
\left(
   \begin{array}{cc}
    \cosh l_{ij}\cosh^2r_j & \cosh r_i\cosh r_j  \\
     \cosh r_i\cosh r_j & \cosh l_{ij}\cosh^2r_i  \\
    \end{array}
 \right),
\end{aligned}
\end{equation*}
which implies that $\mathrm{Hess}_r \widetilde{G}|_{r=\overline{r}}<0$ and hence the function $\widetilde{G}$ is a strictly concave function on a convex neighborhood $U$ of $\overline{r}$.
Set
\begin{equation*}
f(t)=\widetilde{G}((1-t)r+t\overline{r}),\ t\in [0,1].
\end{equation*}
Then $f(t)$ is a $C^1$ smooth concave function of $t\in [0,1]$.
The rest of the proof is similar to that of Theorem \ref{Thm: local rigidity}, we omit it here.
\qed

\end{document}